\documentclass[12pt,leqno]{amsart}
\usepackage{amssymb}
\usepackage{amsmath,amssymb,color}
\numberwithin{equation}{section}

\newcommand{\ep}{\varepsilon}
\newcommand{\la}{\lambda}
\newcommand{\va}{\varphi}
\newcommand{\ppp}{\partial}

\newcommand{\ddd}{\mbox{div}\thinspace}
\newcommand{\rrr}{\mbox{rot}\thinspace}

\newcommand{\weight}{e^{2s\va}}
\newcommand{\sumj}{\sum_{j=1}^3}

\newcommand{\R}{\mathbb{R}}

\newcommand{\N}{\mathbb{N}}

\newcommand{\www}{\widetilde}
\newcommand{\wdel}{\delta_2}
\newcommand{\weps}{\www{\ep}}

\newcommand{\ooo}{\overline}
\newcommand{\OOO}{\Omega}

\newcommand{\dddd}{\Vert d\Vert_{C(\ooo{\OOO})}}

\allowdisplaybreaks

\title{Carleman estimate for the Navier-Stokes equations and
applications}

\author{
$^1$ Oleg Y. Imanuvilov, $^2$ Luca Lorenzi,
$^{3,4,5,6}$ M.~Yamamoto }
\thanks{
$^1$ Department of Mathematics, Colorado State University\\
101 Weber Building, Fort Collins, CO 80523-1874, USA\\
e-mail: {\tt oleg@math.colostate.edu}\\
$^2$ Department of Mathematical, Physical and Computer Sciences, 
University of Parma, Parco Area delle Scienze 53/A, I-43124 Parma
Italy\\
e-mail:{\tt luca.lorenzi@unipr.it}\\
$^3$ Graduate School of Mathematical Sciences, The University
of Tokyo, Komaba, Meguro, Tokyo 153-8914, Japan \\
$^4$ Honorary Member of Academy of Romanian Scientists, 
Ilfov, nr. 3, Bucuresti, Romania \\
$^5$ Correspondence member of Accademia Peloritana dei Pericolanti,\\
Palazzo Universit\`a, Piazza S. Pugliatti 1 98122 Messina Italy \\
$^6$ Peoples' Friendship University of Russia 
(RUDN University) 6 Miklukho-Maklaya St, Moscow, 117198, Russian Federation
e-mail: {\tt myama@ms.u-tokyo.ac.jp}
}

\date{}
\begin{document}
\maketitle

\begin{abstract}
For linearized Navier-Stokes equations, we first derive a Carleman 
estimate with a regular weight function.
Then we apply it to establish conditional stability  for the 
lateral Cauchy problem and finally we prove conditional 
stability estimates for inverse source problem of determining 
a spatially varying divergence-free factor of a source term.
\end{abstract}

\section{Introduction}

Let $T > 0$, and $\OOO \subset \R^3$ be a bounded domain with smooth boundary
$\ppp\OOO$ and $A$, $B$ be sufficiently smooth.
By $a^T$ we denote the transpose of a vector $a$ under consideration.
We set 
$$
\ppp_j = \frac{\ppp}{\ppp x_j}, \quad j=1,2,3, \quad 
\ppp_t = \frac{\ppp}{\ppp t}, \quad 
\nabla = (\ppp_1, \ppp_2, \ppp_3), \quad \nabla_{x,t} = (\nabla, \ppp_t)
$$
and
$$
\ddd v = \sum_{j=1}^3 \ppp_jv_j \quad \mbox{for $v=(v_1,v_2,v_3)^T$}.
$$
 
We consider
$$
\left\{\begin{array}{rl}
& \ppp_tv(x,t) - \Delta v + (A(x,t)\cdot\nabla)v
+ (v\cdot \nabla)B(x,t) + \nabla p = F(x,t),\\        
& \ddd v(x,t) = 0, \quad x \in \OOO, \, 0<t<T. 
\end{array}\right.
                                                         \eqno{(1.1)}
$$

Let $\Gamma \subset \ppp\OOO$ and $\OOO_0 \subset
\OOO$ be given subboundary and subdomain.
\\

In this article, we consider the following two problems.
\\
{\bf Continuation of solution by Cauchy data on $\Gamma \times (0,T)$}.
{\it 
Let $(v,p)$ satisfy (1.1) with $F=0$ and let $I \subset (0,T)$.
Then determine $(v,p)$ in $\OOO_0 \times I$ by $v$ and its derivatives
on $\Gamma \times (0,T)$.
}
\\
\vspace{0.2cm}
{\bf Inverse source problem}.
{\it 
Let $t_0 \in (0,T)$ be given and let $(v,p)$ satisfy (1.1).  
Then determine $F$ 
by $v$ and its derivatives on $\Gamma \times (0,T)$ and 
$v(\cdot,t_0)$.
}
\\

A method based on Carleman estimate by Bukhgeim and Klibanov \cite{BK} is 
applicable for proving stability estimates for these problems.
Since \cite{BK}, there have been many works for parabolic, hyperbolic
and other types of equations.  We refer only to 
Beilina and Klibanov \cite{BeKl}, Bellassoued and Yamamoto \cite{BY},
Imanuvilov and Yamamoto \cite{IY1}, \cite{IY2}, \cite{IY3}, 
Klibanov \cite{Kli}, Klibanov and Timonov \cite{KT}, 
Yamamoto \cite{Ya}.
Here the references are very limited and see also the references in the 
above works.

On the other hand, there are not sufficient works on inverse problems
by the same methodology 
for the Navier-Stokes equations in spite of the importance.
We refer only to 
Imanuvilov and Yamamoto \cite{IYNS}, Bellassoued, Imanuvilov and Yamamoto 
\cite{BIY}, Boulakia \cite{Bou}, 
Choulli, Imanuvilov, Puel and Yamamoto \cite{CIPY},
Fabre \cite{Fab}, Fan, Di Cristo, Jiang and Nakamura \cite{FDJN},
Fan, Jiang and Nakamura \cite{FJN}.  See also Imanuvilov and 
Yamamoto \cite{IY3}, \cite{IY21} 
as for inverse problems for the fluid equations.

Huang, Imanuvilov and Yamamoto \cite{HIY} recently 
simplified the method by \cite{BK} and here we apply it.
We remark that we use the same type of Carleman estimate for the 
Navier-Stokes equations as \cite{BIY} but the derivation 
of a Carleman estimate is simpler than \cite{BIY}.

Thus, in this article, we aim at improving the results in \cite{BIY}, and the 
main achievements are:
\begin{itemize}
\item
a more feasible Carleman estimate (Theorem 1) for the applications.
\item
an improvement of the conditional stability estimates for the continuation of
solutions by Cauchy data.
\item
novel conditional stability estimates in determining force term $F(x,t)$
in view of the divergence-free component. 
\end{itemize}

Now we show the key Carleman estimate, and to this end, we need 
some notations.  Let $t_0 \in (0,T)$ be arbitrarily taken.
We choose $\delta>0$ such that $0<t_0-\delta < t_0+\delta < T$ and we
set 
$$
I = (t_0-\delta, \, t_0+\delta), \quad
Q:= \OOO \times I.
$$
Let $\gamma = (\gamma_1, \gamma_2, \gamma_3) \in (\N \setminus \{0\})^3$ 
and $\vert \gamma\vert = \gamma_1 + \gamma_2 + \gamma_3$.
For $k, \ell \in \N \cup \{0\}$ and $s>0$, we set
$$
H^{k,\ell}(Q) = \{ v \in L^2(Q); \thinspace
\ppp_x^{\gamma}v \in L^2(Q), \thinspace \vert \gamma\vert \le k,
\thinspace \ppp_t^{j}v \in L^2(Q) \quad 0\le j\le \ell\}.
$$

Let $d \in C^2(\ooo{\OOO})$ satisfy 
$$
\vert \nabla d\vert \ne 0 \quad \mbox{on $\ooo{\OOO}$}.
                                          \eqno{(1.2)}
$$
Henceforth we fix a large constant $\la > 0$, and throughout
the article, we omit the $\la$-dependency of the constants. 
For arbitrarily chosen constant $\beta > 0$, we set 
$$
\va(x,t) = e^{\la(d(x) - \beta (t-t_0)^2)}, \quad (x,t) \in Q.
$$
Henceforth by $[u]_k$ with $k=1,2,3$,
we denote the $k$-th component of a vector $u$, and 
$\nabla(\rrr v)$ means the matrix with the elements 
$\ppp_j[\rrr v]_k$, $1\le j,k \le 3$. 

Then 
\\
{\bf Theorem 1 (Carleman estimate for the Navier-Stokes equations).}
\\
{\it 
There exist constants $s_0>0$ and $C>0$ such that 
\begin{align*}
& \int_Q \biggl\{
\frac{1}{s}( \vert \ppp_t\rrr v\vert^2
+ \vert \ppp_tv\vert^2 + \vert \Delta (\rrr v)\vert^2
+ \vert \Delta v\vert^2)\\
+& s(\vert \nabla(\rrr v)\vert^2 + \vert \nabla v\vert^2)
+ s^3(\vert \rrr v\vert^2 + \vert v\vert^2) \biggr\}
\weight dxdt \\
\le& C\int_Q \vert \rrr F\vert^2 \weight dxdt 
+ Cs^3\int_{\ppp\OOO\times I} 
\sum_{j=0}^1 (\vert \nabla_{x,t}^j (\rrr v) \vert^2 
+ \vert \nabla^jv\vert^2) \weight dSdt\\
+ &Cs^3 \sum_{\kappa=0}^1 \int_{\OOO} \sum_{j=0}^1 
\vert \nabla^j (\rrr v(x,t_0+(-1)^{\kappa}\delta))
\vert^2 e^{2s\va(x,t_0+\delta)} dx
\end{align*}
for all $s > s_0$, $p\in L^2(I;H^1(\OOO))$ and $v, \rrr v\in H^{2,1}(Q)$.
}
\\

We note that $\va(x,t_0+\delta) = \va(x,t_0-\delta)$, $x\in \OOO$.
\\
{\bf Remark 1.}
\\
We can estimate also $\frac{1}{s}\vert \nabla p\vert^2 \weight$ on the 
left-hand side if we replace the term 
$\int_Q \vert \rrr F\vert^2 dxdt$ by 
$\int_Q (\vert F\vert^2 + \vert \rrr F\vert^2 dxdt$. 
\\

Here we choose a regular weight function $e^{\la(d(x) - \beta (t-t_0)^2)}$.
A similar Carleman estimate with the same type of weight function for 
the Navier-Stokes equations is derived in \cite{BIY}, but 
the derivation here is simpler.
On the other hand, we can choose other type of weight function with 
singularity at $t=0, T$, which was created by Imanuvilov \cite{Ima} for 
a single parabolic equation.  For linearized Navier-Stokes equations,
we refer to Choulli, Imanuvilov, Puel and Yamamoto \cite{CIPY}, 
Fan, Di Cristo, Jiang and Nakamura \cite{FDJN}, Fan, Jiang and Nakamura 
\cite{FJN}, Fern\'andez-Cara, Guerrero, Imanuvilov and Puel \cite{FGIP}.
\\

Next we formulate the main result for the continuation of solution.
For small $\delta>0$ satisfying $0<t_0-\delta < t_0+\delta$, we choose 
$0 < \delta_0 < \delta<T$. 
We consider 
$$
\left\{\begin{array}{rl}
& \ppp_tv(x,t) - \Delta v + (A(x,t)\cdot\nabla)v + (v\cdot \nabla)B
+ \nabla p = 0,  \\
& \ddd v = 0, \quad (x,t) \in Q:= \OOO \times I. 
\end{array}\right.
                                                         \eqno{(1.3)}
$$
Then we are requested to estimate $(v,p)$ in some subdomain by 
data of $v$ on $\Gamma \times (0,T)$, where $\Gamma$ is an arbitrary 
subboundary of $\ppp\OOO$.

For an arbitrarily chosen subboundary $\Gamma \subset \ppp\OOO$, we choose 
an arbitrary subdomain $\OOO_0 \subset\OOO$ such that
$$
\ooo{\OOO_0} \subset \OOO \cup \Gamma, \quad
\mbox{$\ooo{\OOO_0} \cap \ppp\OOO$ is included in Int $(\Gamma)$: the interior
of $\Gamma$.}
$$

We are ready to state the main result for the continuation.
\\
{\bf Theorem 2 (continuation of solution).}\\
{\it
Let $(v,p) \in H^{2,1}(Q) \times L^2(0,T;H^2(\OOO))$ satisfy (1.3) 
with $\rrr v\in H^{2,1}(Q)$ and 
$$
\sum_{j=0}^1 (\Vert \nabla_{x,t}^j(\rrr v)\Vert_{L^2(\ppp\OOO\times (0,T))} 
+ \Vert \nabla^jv\Vert_{L^2(\ppp\OOO\times (0,T))})
+ \Vert \rrr v\Vert_{L^{\infty}(0,T; H^1(\OOO))} \le M            \eqno{(1.4)}
$$
with arbitrarily fixed constant $M>0$.
Then for any $\ep>0$,
there exist constants $C>0$ and $\theta\in (0,1)$, which depend on 
$M, \ep, \Gamma, \OOO_0$ such that 
$$
\sum_{\ell=0}^1 \Vert \rrr^{\ell}v\Vert_{H^{1,1}(\OOO_0\times (\ep, \,T-\ep))}
+ \Vert \Delta (\rrr^{\ell}v)\Vert_{L^2(\OOO_0\times (\ep,\,T-\ep))}
+ \Vert \nabla p\Vert_{L^2(\OOO_0\times (\ep, T-\ep))}\\
\le C( D + D^{\theta}),
$$
where
$$
D:= \sum_{j=0}^1 (\Vert \nabla_{x,t}^j(\rrr v)\Vert
_{L^2(\Gamma\times (0,T))} 
+ \Vert \nabla^jv\Vert_{L^2(\Gamma\times (0,T))}).
$$
}
\\

This theorem asserts stability in determining a solution $(v,p)$ 
in some subdomain by boundary data of $v$ on $\Gamma \times (0,T)$
under a priori bound (1.4).

A similar conditional stability for the continuation was proved in 
\cite{BIY}, but our proof is simplified by adapting the proof of 
Proposition 2 in Huang, Imanuvilov and Yamamoto \cite{HIY} which is concerned 
with a single parabolic equation.
In \cite{HIY}, differently from the conventional arguments by 
Carleman estimates, we do not use any cut-off functions $\chi(x,t)$ 
which is compactly supported, and do not need to consider 
$\chi(x,t)v(x,t)$, but directly $v(x,t)$.  In the case of the 
Navier-Stokes equations, the cut-off destroys the original 
structure of the equations, that is, $\ddd \chi v = 0$ does not hold.
Consequently, if we apply the cut-off, then we have to prove
a Carleman estimate for the Navier-Stokes equations without 
$\ddd v=0$, which causes a complexity (see \cite{BIY}).  
Since we do not here require any cut-off at all,
the proofs of Theorems 1 - 3 are simplified.

Finally we state the main result for the inverse source problem.
In (1.1), we assume
$$
\ddd F(x,t_0) = 0, \quad x\in \OOO.            \eqno{(1.5)}
$$
We introduce two kinds of conditions:
$$
\vert \ppp_t^k\rrr F(x,t)\vert \le C\vert \rrr F(x,t_0)\vert, 
\quad (x,t) \in Q, \, k=0,1.                                 \eqno{(1.6)}
$$
$$
\vert \ppp_t^k\rrr F(x,t)\vert \le C(\vert \nabla F(x,t_0)\vert
+ \vert F(x,t_0)\vert), \quad (x,t) \in Q, \, k=0,1,2.
                                                                 \eqno{(1.7)}
$$
We notice that (1.6) implies (1.7), that is, condition (1.7) is more 
generous than (1.7).
\\

We are ready to state different conditional stability estimates 
according to 
(1.6) and (1.7).
\vspace{0.2cm}
\\
{\bf Theorem 3 (inverse source problem).}\\
{\it 
Let $\Gamma \subset \ppp\OOO$ be arbitrary and 
$\OOO_0 \subset \OOO$ be a subdomain satisfying $\ooo{\OOO_0}
\subset \OOO \cup \Gamma$ and $\ooo{\OOO_0} \cap \ppp\OOO
\subset \mbox{Int}\, (\Gamma)$. 
\\
(i) We set 
$$
\mathcal{F}_1:= \{ F \in L^2(Q);\, \rrr F, \ddd F \in H^1(I;L^2(\OOO)),
\, \mbox{(1.6) is satisfied} \}.
$$
We set 
$$
D_1:= \sum_{j,k=0}^1 (\Vert \nabla_{x,t}^j\ppp_t^k(\rrr v)\Vert
_{L^2(\Gamma\times I)} + \Vert \nabla^j\ppp_t^kv\Vert
_{L^2(\Gamma\times I)}) + \Vert v(\cdot,t_0)\Vert_{H^3(\OOO)}.
                                                       \eqno{(1.8)}
$$
We assume that a solution $(v,p) \in H^{2,1}(Q) \times H^{1,0}(Q)$ to (1.1)
satisfies $\ppp_t\rrr v\in H^{2,1}(Q)$, $\ppp_tp \in H^{1,0}(Q)$ and 
$$
\sum_{j,k=0}^1 (\Vert \nabla_{x,t}^j\ppp_t^k(\rrr v)\Vert
_{L^2(\ppp\OOO\times I)} 
+ \Vert \nabla^j\ppp_t^kv\Vert_{L^2(\ppp\OOO\times I)})
$$
$$
+ \sum_{\kappa=0}^1 
\sum_{k=0}^1 \Vert \ppp_t^k\rrr (\cdot,\,t_0+(-1)^{\kappa}\delta)
\Vert_{H^1(\OOO)}
\le M_1                                                \eqno{(1.9)}
$$
with arbitrarily given constant $M_1>0$. 
Then there exist constants $C>0$ and $\theta \in (0,1)$, dependent 
on $\OOO_0, \Gamma, M_1$, such that 
$$
\Vert \rrr F(\cdot,t_0)\Vert_{L^2(\OOO_0)} 
\le C(D_1^{\theta} + D_1)
$$
for each $F \in \mathcal{F}_1$.
\\
(ii) We set 
\begin{align*}
&\mathcal{F}_2:= \{ F \in H^2(I;H^2(\OOO));\, 
\Vert F(\cdot,t_0)\Vert_{H^2(\OOO)} \le m_1, \quad
F(\cdot,t_0) = \vert \nabla F(\cdot,t_0)\vert = 0
\quad \mbox{on $\Gamma$}, \\
& \quad \mbox{(1.5) and (1.7) are satisfied} \}
\end{align*}
with arbitrarily chosen constant $m_1>0$.
We further set
$$
D_2:= \sum_{j=0}^1 \sum_{k=0}^2
(\Vert \nabla_{x,t}^j\ppp_t^k(\rrr v)\Vert_{L^2(\Gamma\times I)} 
+ \Vert \nabla^j\ppp_t^kv\Vert
_{L^2(\Gamma\times I)}) + \Vert v(\cdot,t_0)\Vert_{H^4(\OOO)}.
$$
We assume that a solution $(v,p) \in H^{2,1}(Q) \times H^{1,0}(Q)$ to 
(1.1) satisfies $\ppp_t^k\rrr v \in H^{2,1}(Q)$, $\ppp_t^kp \in H^{1,0}(Q)$ 
with $k=0,1,2$ and 
\begin{align*}
&\sum_{j=0}^1 \sum_{k=0}^2 
(\Vert \nabla_{x,t}^j\ppp_t^k\rrr v\Vert_{L^2(\ppp\OOO\times I)}
+ \Vert \nabla^j \ppp_t^kv\Vert_{L^2(\ppp\OOO\times I)})\\
+& \sum_{\kappa=0}^1\sum_{k=0}^2 
\Vert \ppp_t^k\rrr v(\cdot,\,t_0+(-1)^{\kappa}\delta)\Vert
_{H^1(\OOO)} \le M_2
\end{align*}
with arbitrarily given constant $M_2>0$. 
Then there exist constants $C>0$ and $\theta \in (0,1)$, dependent 
on $\OOO_0, \Gamma, m_2, M_2$, such that 
$$
\Vert F(\cdot,t_0)\Vert_{H^1(\OOO_0)} \le C(D_2^{\theta} + D_2)
$$
for each $F \in \mathcal{F}_2$.
}
\\

In Theorem 3 (i) we aim at the determination of rot $F(x,t_0)$ for 
$x\in \OOO$, not $F$ itself, and we do not need to assume (1.5). 
Theorem 3 asserts two kinds of conditional stability according to 
the admissible sets $\mathcal{F}_1$ and $\mathcal{F}_2$.
Theorem 3 (i) determines only the rotation component of 
$F(x,t_0)$ and so it is not necessary to assume (1.5).

In Theorem 3, not assuming the boundary condition on the whole
boundary of $\OOO$, we establish stability in determining $F$ in a subdomain 
of $\OOO$.  On the other hand, the stability over $\OOO$ is proved with the 
boundary condition on the whole $\ppp\OOO \times (0,T)$ in 
Choulli, Imanuvilov, Puel and 
Yamamoto \cite{CIPY} for the Navier-Stokes equations, 
and Fan, Jiang and Nakamura \cite{FJN} for the Boussinesq equations.
\\

We will test conditions (1.5) - (1.7).
First, as the following example shows,
the condition (1.5) is essential for the uniqueness in 
determining $F(x,t_0)$ in Thoerem 3 (ii), where no data of $p$ are 
observed.  
\\
{\bf Obstruction to the uniqueness in determining $F$:}
\\
We consider a simple case
$$
\left\{\begin{array}{rl}
& \ppp_tv(x,t) - \Delta v + \nabla p = f(x),  \\
& \ddd v = 0, \quad (x,t) \in Q:= \OOO \times I,\\
& v(x,t_0) = 0 \quad\mbox{for $x\in \OOO$}, \quad
\mbox{supp $v \subset \OOO \times (0,T)$}. 
\end{array}\right.
                                                         \eqno{(1.10)}
$$
Here $f$ is an $\R^3$-valued smooth function.
It is trivial that $(v,p) = (0,0)$ satisfies (1.10) with $f=0$.
Let $\psi\in C^{\infty}_0(\OOO)$.  Then $(v,p) := (0, \, \psi)$
satisfies (1.10) with $f:= \nabla \psi$.  In other words, by the 
appearance of the pressure field $p$ in the Navier-Stokes equations, 
there is no possibility for the uniqueness in 
determining the component of $f$ given by a
scalar potential.
\\

Next we examine the conditions (1.5) - (1.7) by considering examples.
\\
{\bf Example 1.}
We see that (1.5) is replaced in terms of a vector potential:
$$
F(x,t) := \rrr q(x,t), \quad (x,t) \in Q,
$$
where $q = q(x,t)$ is smooth.  Then (1.5) is automatically satisfied.
Moreover, we see that (1.7) is equivalent to
$$
\vert \rrr \rrr (\ppp_t^kq(x,t))\vert
\le C(\vert \nabla \rrr q(x,t_0)\vert + \vert \rrr q(x,t)\vert),\quad
(x,t) \in Q,\, k=0,1,2.
$$
\\
{\bf Example 2.}\\
Let 
$$
F(x,t) = r(t)f(x),  \quad x\in \OOO, \, 0<t<T,  
\quad r(t_0) \ne 0, \quad \ddd f = 0 \quad \mbox{in $\OOO$},
$$
where real-valued $r\in C^2[0,T]$ and $f=(f_1,f_2,f_3)^T$ are smooth.
Then (1.5) and (1.6) are satisfied.  Indeed we can directly verify 
$\ddd F(x,t_0) = r(t_0)\ddd f(x) = 0$ for $x\in \OOO$, which is (1.5).
Moreover, we have
$$
\vert \ppp_t^k\rrr F(x,t)\vert 
= \left\vert \frac{d^kr}{dt^k}(t) \rrr f(x)\right\vert 
\le C\vert \rrr f(x)\vert, \quad k=0,1
$$
and
$$
\vert \rrr F(x,t_0)\vert = \vert r(t_0)\vert \vert \rrr f(x)\vert, \quad
(x,t) \in Q,
$$
which is (1.6).  
\\
{\bf Example 3.}\\
For more general $F$ unlike Example 2, it is not always easy to verify 
(1.6), but in the following case the condition (1.7) can be readily 
verified. 
Let 
$$
F(x,t) = R(x,t)f(x),
$$
where $R(x,t) = (r_{jk}(x,t))_{1\le j,k\le 3}$ and each component $r_{jk}$ is 
smooth in $(x, t)$, and $f(x) = (f_1(x), f_2(x), f_3(x))^T$.
Then (1.7) is satisfied if 
$$
\vert \mbox{det}\, R(x,t_0)\vert > 0 \quad \mbox{for $x\in \ooo{\OOO}$}.
$$

Indeed, we have
$$
\vert \ppp_t^k \rrr F(x,t) \vert \le C(\vert \nabla f(x)\vert 
+ \vert f(x)\vert),\quad (x,t) \in \ooo{Q}, \quad k=0,1,2.
$$
Therefore it is sufficient to verify 
$$
\vert \nabla f(x)\vert + \vert f(x)\vert
\le C(\vert \nabla F(x,t_0)\vert + \vert F(x,t_0)\vert),
\quad x \in \OOO.                     \eqno{(1.11)}
$$
{\bf Verification of (1.11).}
By $\vert \mbox{det}\, R(x,t_0)\vert \ne 0$ for $x \in \ooo{\OOO}$, we have 
$$
\vert f(x)\vert \le C\vert F(x,t_0)\vert, \quad x\in \ooo{\OOO}.
$$
Next, setting $F = (F_1,F_2,F_3)^T$, we have
$$
\ppp_{\ell}F_k(x,t_0) = \sumj r_{kj}(x,t_0)\ppp_{\ell}f_j(x)
+ \sumj (\ppp_{\ell}r_{kj})(x,t_0)f_j(x),
$$
and so 
$$
\sumj r_{kj}(x,t_0) \ppp_{\ell}f_j(x) = \ppp_{\ell}F_k(x,t_0)
- \sumj (\ppp_{\ell}r_{kj}(x,t_0))f_j(x),
$$
that is,
$$
R(x,t_0)\ppp_{\ell}f(x) = \ppp_{\ell}F(x,t_0)
- (\ppp_{\ell}R(x,t_0))f(x).
$$
Therefore, by $\vert \mbox{det}\, R(x,t_0)\vert \ne 0$ for $x\in \ooo{\OOO}$, 
we see
$$
\vert \ppp_{\ell}f(x) \vert \le C(\vert \ppp_{\ell}F(x,t_0)\vert 
+ \vert F(x,t_0)\vert), \quad x \in \ooo{\OOO}, \quad 1\le \ell\le 3.
$$
Thus (1.11) is verified.
\\
\vspace{0.2cm}

Here for the form $F(x,t) = R(x,t)f(x)$ with $3\times 3$ matrix 
$R(x,t)$, we do not know a convenient sufficient condition for (1.6),
but condition (1.7)
can be more directly verified.
Therefore, Theorem 3 (ii) is more feasible in applications.
\\
{\bf Remark 2.}\\
A special case $F(x,t) = f(x)r(x,t)$, where 
$f$ is real-valued and $r$ is $\R^3$-valued, is trivial by 
(1.5).  Because (1.5) means that $f$ satisfies a first-order 
transport equation $r(x,t_0)\cdot \nabla f(x) + f(x)\ddd r(x, t_0)=0$ for
$x \in \OOO$, and so $f$ can be uniquely determined by boundary data of 
$f$ under a certain condition on $r(x,t_0)$.
\\
{\bf Remark 3.}
\\
Theorem 3 (i) estimates only $\rrr F(x,t_0)$.  However, since
$\ddd F(x,t_0) = 0$ for $x\in \OOO$, we have 
$$
-\Delta F(x,t_0) = \rrr \rrr F(x,t_0) - \nabla \ddd F(x,t_0)
= \rrr \rrr F(x,t_0), \quad x\in \OOO.
$$
If we are given boundary data $F(x,t_0)$ with $x \in \ppp\OOO$ in 
a suitable space, then we can further estimate by means of 
$H^{-1}$-estimate of the solution to $-\Delta F(x,t_0)=h$ under 
condition (1.6).  Here we do not pursue more.  In Section 5.4, we will 
discuss this issue by means of a Carleman estimate in $H^{-1}$-space
under a weaker assumption than (1.7). 

This article is composed of six sections.
In Sections 2-4, we prove Theorems 1 -3 respectively.
Section 5 is devoted to concluding remarks and in Section 6, we 
derive Lemma 2 which is used for the proof of Theorem 1.

\section{Proof of Theorem 1}

Thanks to the large parameter $s>0$, it suffices to prove 
Theorem 1 for $B=0$.

We set 
$$
z := \rrr v.
$$
Then 
$$
\left\{\begin{array}{rl}
& \ppp_tz(x,t) - \Delta z + (A(x,t)\cdot\nabla)z
= \rrr F - \sumj \nabla A_j \times \ppp_jv, \\
& \Delta v = -\rrr z, \quad x \in \OOO, \, 0<t<T. 
\end{array}\right.
                                                         \eqno{(2.1)}
$$
\\
The system (2.1) provides a decomposition of the Navier-Stokes equations into 
a parabolic equation in $z:= \rrr v$ and an elliptic equation in $v$, where 
$p$ is eliminated.  

First we show a Carleman estimate for a simple parabolic equation with constant
coefficients:
$$
\ppp_tu - \Delta u = G \quad \mbox{in $Q$}.
$$

Henceforth $C>0$ denotes generic constants which are independent of 
$s>0$.
\\
\vspace{0.2cm}
{\bf Lemma 1.}\\
{\it
There exist constants $C>0$ and $s_0>0$ such that 
\begin{align*}
& \int_Q \left\{
\frac{1}{s}( \vert \ppp_t u\vert^2 + \vert \Delta u\vert^2)
+ s\vert \nabla u\vert^2 + s^3\vert u\vert^2 \right\}
\weight dxdt \\
\le& C\int_Q \vert G\vert^2 \weight dxdt 
+ Cs^3\int_{\ppp\OOO \times I} (\vert \nabla_{x,t}u\vert^2
+ \vert u\vert^2) \weight dSdt\\
+& Cs^3 \sum_{\kappa=0}^1 \int_{\OOO} 
(\vert \nabla u(x,t_0 + (-1)^{\kappa}\delta)\vert^2
+ \vert u(x,t_0 + (-1)^{\kappa} \delta)\vert^2) e^{2s\va(x,t_0+\delta)} dx
\end{align*}
for all $s > s_0$ and all $u\in H^{2,1}(Q)$.
}
\\

Lemma 1 is a standard Carleman estimate for a single parabolic equation,
and a direct proof is found 
for example in Bellassoued and Yamamoto \cite{BY} (Lemma 7.1),
Yamamoto \cite{Ya} (Theorems 3.1 and 3.2).
\\

Next we show two Carleman estimates for the Laplace operator.
\\
{\bf Lemma 2.}\\
{\it
Let $r \in L^2(I;H^2(\OOO))$ satisfy 
$$
-\Delta r(x,t) = g(x,t) \quad \mbox{in $Q$.}
$$
Then
\begin{align*}
& \int_Q (s\vert \nabla r(x,t)\vert^2 + s^3\vert r(x,t)\vert^2)e^{2s\va(x,t)}
dxdt \\
\le & C\int_Q \vert g(x,t)\vert^2 e^{2s\va(x,t)} dxdt 
+ Cs^3\int_{\ppp\OOO\times I} (\vert \nabla r\vert^2 + \vert r\vert^2)
e^{2s\va(x,t)} dSdt
\end{align*}
for all $s > s_0$.
}
\\
{\bf Lemma 3.}\\
{\it
Let $w \in H^2(\OOO)$ satisfy
$$
-\Delta w(x) = h(x), \quad x\in \OOO.
$$
Then 
\begin{align*}
& \int_{\OOO} (s\vert \nabla w(x)\vert^2 + s^3\vert w(x)\vert^2)
e^{2s\va(x,t_0)} dx\\
\le &C\int_{\OOO} \vert h(x)\vert^2 e^{2s\va(x,t_0)} dx 
+ Cs^3\int_{\ppp\OOO} (\vert \nabla w(x)\vert^2 + \vert w(x)\vert^2)
e^{2s\va(x,t_0)} dS
\end{align*}
for all $s>s_1$.
}
\\

The proof of Lemma 3 is standard and we can prove directly by integration by 
parts (e.g., Lemma 7.1 in \cite{BY}), and Lemma 3 is used also for the 
proof of Theorem 3 in Section 4.

The proof of Lemma 2 is done by integrating the standard elliptic Carleman 
estimate in Lemma 3 for $-\Delta$ over the time interval $I$.  
For completeness, the derivation of 
Lemma 2 from Lemma 3 is provided in Appendix.
\\

Now we complete the proof of Theorem 1.
We apply Lemmata 1 and 2 to (2.1):
\begin{align*}
&\int_Q \left(\frac{1}{s} (\vert \ppp_tz\vert^2 + \vert \Delta z\vert^2)
+ s\vert \nabla z\vert^2 + s^3\vert z\vert^2 
\right\}\weight dxdt\\
\le & C\int_Q \left\vert \rrr F - \sumj A_j\times (\ppp_jv)
\right\vert^2 \weight dxdt + CJ_1
\end{align*}
$$
\le C\int_Q \vert \rrr F \vert^2 \weight dxdt
+ C\int_Q (\vert \nabla v\vert^2 + \vert v\vert^2) \weight dxdt
+ CJ_1                               \eqno{(2.2)}
$$
and
$$
\int_Q (s\vert \nabla v\vert^2 + s^3\vert v\vert^2) \weight dxdt 
\le C\int_Q \vert \rrr z \vert^2 \weight dxdt + CJ_2.
                                           \eqno{(2.3)}
$$
Here and henceforth we set 
\begin{align*}
&J_1 := s^3\int_{\ppp\OOO\times I}
(\vert \nabla_{x,t}(\rrr v)\vert^2 + \vert \rrr v\vert^2)\weight dSdt\\
+& s^3\sum_{\kappa=0}^1\int_{\OOO} 
(\vert \nabla(\rrr v(x,t_0+(-1)^{\kappa}\delta))\vert^2
+ \vert \rrr v(x,t_0+(-1)^{\kappa}\delta)\vert^2)e^{2s\va(x,t_0+\delta)} dx
\end{align*}
and
$$
J_2 := s^3\int_{\ppp\OOO\times I}
(\vert \nabla v\vert^2 + \vert v\vert^2)\weight dSdt.
$$
Hence, substituting (2.3) into the second term on the right-hand side
of (2.1):
\begin{align*}
& \int_Q \left\{
\frac{1}{s} (\vert \ppp_tz\vert^2 + \vert \Delta z\vert^2)
+ s\vert \nabla z\vert^2 + s^3\vert z\vert^2 
\right\}\weight dxdt \\
\le & C\int_Q \vert \rrr F \vert^2 \weight dxdt
+ C\int_Q \vert \rrr z\vert^2 \weight dxdt + C(J_1+J_2).
\end{align*}
Absorbing the second term on the right-hand side into the left-hand side,
we complete the proof of Theorem 1.

\section{Proof of Theorem 2}

Theorem 1 holds with general $d(x)$ satisfying only (1.2) 
in the weight function, but
in order to apply for the proof of Theorems 2 and 3, we have to choose 
$d(x)$ closely related to the geometry of $\Gamma$ as follows.
First we construct some domain $\OOO_1$.
For $\Gamma \subset \ppp\OOO$, we choose a bounded domain $\OOO_1$ 
with smooth boundary such that
$$
\OOO \subsetneqq \OOO_1, \quad \ooo{\Gamma} = \ooo{\ppp\OOO\cap\OOO_1}, 
\quad \ppp\OOO\setminus\Gamma \subset \ppp\OOO_1.
$$
In particular, $\OOO_1\setminus\ooo{\OOO}$ contains some non-empty 
open subset. 
We note that $\OOO_1$ can be constructed as the interior of a union  
of $\ooo{\OOO}$ and the closure of a non-empty domain 
$\widehat{\OOO}$ satisfying $\widehat{\OOO} \subset \R^n 
\setminus \ooo{\OOO}$ and $\ppp\widehat{\OOO} \cap \ppp\OOO = \Gamma$.

We choose a domain $\omega$ such that 
$\ooo{\omega} \subset \OOO_1 \setminus \ooo{\OOO}$.
Then, by Imanuvilov \cite{Ima}, we can find $d\in C^2(\ooo{\OOO_1})$ such that
$$
d>0 \quad \mbox{in $\OOO_1$}, \quad |\nabla d| > 0 \quad \mbox{on
$\ooo{\OOO_1\setminus \omega}$}, \quad d=0 \quad 
\mbox{on $\ppp\OOO_1$}.                               \eqno{(3.1)}
$$
In particular, by $\ooo{\OOO_0} \subset \OOO_1$, we see
$$
d > 0 \quad \mbox{on $\ooo{\OOO_0}$}, \quad 
d=0 \quad \mbox{on $\ppp\OOO \setminus \Gamma$}.    \eqno{(3.2)}
$$
We recall that we choose a domain $\OOO_0 \subset \OOO$ satisfying 
$\ooo{\ppp{\OOO_0} \cap \ppp\OOO} \subset \mbox{Int}(\Gamma)$ and 
$\ooo{\OOO_0} \subset \OOO \cup \Gamma$.
\\

Henceforth we use this $d$ and a sufficiently large constants $\la > 0$.
Later by (3.10), we choose a constant $\beta > 0$.

For given $t_0>0$ and $\delta_2 > 0$ satisfying $0<t_0-\delta_2
<t_0 + \wdel < T$.  We apply the Carleman estimate Theorem 1 in 
$\OOO \times (t_0-\wdel, t_0+\wdel)$ to obtain
\begin{align*}
& \int_{\OOO \times (t_0-\delta_2,t_0+\delta_2)}
\left( \sum_{j,\ell=0}^{1} \vert \nabla_{x,t}^j(\rrr^{\ell}v)\vert^2
+ \sum_{\ell=0}^1 \vert \Delta (\rrr^{\ell} v)\vert^2\right) \weight dxdt\\
\le& Cs^4\int_{\ppp\OOO\times (t_0-\delta_2,t_0+\delta_2)} 
\sum_{j=0}^1 (\vert \nabla_{x,t}^j(\rrr v) \vert^2 
+ \vert \nabla^j v\vert^2) \weight dSdt
\end{align*}
$$
+ Cs^4\sum_{\kappa=0}^1\int_{\OOO} \sum_{j=0}^1 \vert \nabla^j 
(\rrr v(x,t_0+(-1)^{\kappa}\delta_2))
\vert^2 e^{2s\va(x,t_0+\delta_2)} dx           \eqno{(3.3)}
$$
for all $s > s_0$.

Here we note that the constants $C>0$ and $s_0>0$ are independent of 
$t_0$ because the Carleman estimate is invariant by the translation in time
provided that the translated time interval is included in $(0,T)$.

We set 
$$
d_0 = \min_{x\in \ooo{\OOO_0}} d(x), \quad
d_1 = \max_{x\in \ooo{\OOO}} d(x).
$$

Let $\www{\ep}>0$ be given such that $0 < \www{\ep} < \delta_2$. 
Then, by (3.2), we have
$$
\max_{x\in \ooo{\ppp\OOO\setminus \Gamma},t_0-\wdel\le t \le t_0+\wdel}\va(x,t)
\le 1,
\quad \max_{x\in \ooo{\OOO}} \va(x,t_0-\wdel)
= \max_{x\in \ooo{\OOO}} \va(x,t_0+\wdel)
\le e^{\lambda(d_1^2-\beta \wdel^2)}                   \eqno{(3.4)}
$$
and
$$
\min_{x\in \ooo{\OOO_0}, t_0-\weps\le t\le t_0+\weps} \va(x,t)
\ge e^{\lambda(d_0^2-\beta \weps^2)}.                          \eqno{(3.5)}
$$
For concise descriptions, we set 
$$
N(v;t_0-\www{\ep},t_0+\www{\ep}):= \sum_{\ell=0}^{1} (\Vert \rrr^{\ell}v\Vert^2
_{H^{1,1}(\OOO_0\times (t_0-\www{\ep}, \,t_0+\www{\ep}))}
+ \Vert \Delta (\rrr^{\ell}v)\Vert^2
_{L^2(\OOO_0\times (t_0-\www{\ep}, t_0+\www{\ep})}).
$$ 

Therefore, shrinking the integral domain $\OOO \times 
(t_0-\wdel,t_0+\wdel)$ to $\OOO_0 \times (t_0-\weps,t_0+\weps)$ 
on the left-hand side of (3.3), applying the trace theorem and 
the embedding $H^1(0,T;H^1(\OOO)) \subset C([0,T]; H^1(\OOO))$, 
in terms of (1.4), (3.4) and (3.5), we obtain
$$
\exp(2se^{\lambda(d_0^2-\beta \weps^2)})N(v;t_0-\www{\ep},t_0+\www{\ep})
                                \eqno{(3.6)}
$$
\begin{align*}
\le &Cs^4\left(\int_{\Gamma \times (t_0-\wdel,t_0+\wdel)}
+ \int_{(\ppp\OOO \setminus \Gamma) \times (t_0-\wdel,t_0+\wdel)}\right)
\sum_{j=0}^1 (\vert \nabla_{x,t}^j (\rrr v)\vert^2
+ \vert \nabla^jv\vert^2) e^{2s\va} dSdt                             \\
+& Cs^4\sum_{\kappa=0}^1
\int_{\OOO} \sum_{j=0}^1 \vert \nabla^j\rrr(x, t_0+(-1)^{\kappa}\delta_2)
\vert^2
e^{2s\va(x,t_0+\wdel)} dx\\
\le& Cs^4e^{Cs}D^2 + Cs^4e^{2s}M^2
+ Cs^4M^2\exp(2se^{\lambda(d_1^2-\beta \wdel^2)})
\end{align*}
for all $s\ge s_0$.  

Now we will make a specific choice of $\beta > 0$, $\www{\ep}>0$ and 
$\wdel > 0$ such that $0 < \www{\ep} < \delta_2$.  
First fix sufficiently large $N>1$ such that 
$$
N-1 > \frac{d_1^2 - d_0^2}{d_0^2}, \quad 
\mbox{that is,} \quad \frac{d_1^2}{d_0^2} < N.              \eqno{(3.7)}
$$
For given $\ep>0$ in the statement of the theorem, we set  
$$
\weps := \frac{\ep}{N-1}, \quad 
\wdel := N\weps = \frac{N\ep}{N-1} > \www{\ep}.           \eqno{(3.8)}
$$
Since it suffices to consider the case where $\ep>0$ is sufficiently small,
we can assume that $\delta_2>0$ is small, in particular,
$0<\delta_2 < T-\delta_2 < T$.  Then, we can prove  
$$
\frac{d_1^2 - d_0^2}{\wdel^2 - \weps^2} 
< \frac{d_0^2}{\weps^2}.                               \eqno{(3.9)}
$$
Indeed, by (3.8) and $N>1$ we see 
$$
\frac{\delta_2^2}{\www{\ep}^2} = N^2 > N.
$$
Using (3.7), we have
$$
\frac{d_1^2}{d_0^2} < N < \frac{\delta_2^2}{\www{\ep}^2}, \quad
\mbox{that is,}\quad 
\frac{d_1^2-d_0^2}{d_0^2} < \frac{\delta_2^2-\www{\ep}^2}{\www{\ep}^2}, 
$$
which yields (3.9) by multiplying with $\frac{d_0^2}
{\delta_2^2-\www{\ep}^2} > 0$.
\\

Therefore, we can choose $\beta > 0$ such that 
$$
\frac{d_1^2 - d_0^2}{\wdel^2 - \www{\ep}^2} < \beta 
< \frac{d_0^2}{\www{\ep}^2},                         \eqno{(3.10)}
$$
which implies
$$
d_0^2 - \beta \www{\ep}^2 > 0, \quad 
d_0^2 - \beta \www{\ep}^2 > d_1^2 - \beta \delta_2^2,
$$
that is,
$$
\mu_1:= e^{\lambda(d_0^2-\beta \weps^2)}
> \mu_2:= \max \{1, \, e^{\lambda(d_1^2-\beta \wdel^2)}\}.
$$

We arbitrarily choose $t_0 \in (\delta_2, T-\delta_2)$.
We notice that $(t_0-\delta_2, t_0+\delta_2) \subset (0,T)$.
Hence (3.6) yields 
$$
N(v;t_0-\www{\ep}, t_0+\www{\ep}) 
\le Cs^4M^2e^{-2s\mu_0} + Cs^4e^{Cs}D^2
$$
for all $s\ge s_0$.  Here we note 
$$
\mu_0:= \mu_1 - \mu_2 > 0.
$$
Replacing $s$ and $C$ by $s+s_0$ and $e^{Cs_0}$ respectively, we see
$$
N(v;t_0-\www{\ep}, t_0+\www{\ep}) 
\le Cs^4M^2e^{-2s\mu_0} + Cs^4e^{Cs}D^2
$$
for all $s\ge 0$.  By $\mu_0>0$, we easily verify
$\sup_{s>0} s^4e^{-s\mu_0} < \infty$ and 
$s^4e^{Cs} \le C_1e^{C_1s}$ with large $C_1>0$, 
we obtain
$$
N(v;t_0-\www{\ep}, t_0+\www{\ep}) \le CM^2e^{-s\mu_0} + Ce^{Cs}D^2
$$
for all $s\ge 0$.  Here we note that $C>0$ is a generic constant 
which is independent of $s>0$. 
We make the right-hand side small by choosing $s\ge 0$ suitably.
\\
{\bf Case 1: $D<M$.}
\\
We choose $s > 0$ such that 
$$
M^2e^{-s\mu_0} = D^2e^{Cs}, 
$$
that is,
$$
s = \frac{2}{C+\mu_0}\log \frac{M}{D} > 0.
$$
Then 
$$
N(v;t_0-\www{\ep}, t_0+\www{\ep}) \le 2CM^{\frac{2C}{C+\mu_0}}
D^{\frac{2\mu_0}{C+\mu_0}}.
$$
\\
{\bf Case 2: $D \ge M$.}
\\
Then we can directly 
$$
N(v;t_0-\www{\ep}, t_0+\www{\ep}) \le CD^2(e^{-s\mu_0} + e^{Cs}).
$$
Hence setting $s=1$ for example, we see 
$N(v;t_0-\www{\ep}, t_0+\www{\ep}) \le CD^2e^{C}$.
Thus in both cases, we obtain
$$
N(v;t_0-\www{\ep}, t_0+\www{\ep}) \le C(M)(D^2 + D^{2\theta}), 
                                           \eqno{(3.11)}
$$
where $\theta = \frac{\mu_0}{C+\mu_0} \in (0,1)$.
Here the constants $C$ and $\theta$ are dependent on $\www{\ep},
\delta_2 > 0$, but independent of $t_0$.  Varying $t_0$ over 
$(\delta_2, \, T-\delta_2)$, we reach 
$$
N(v;\delta_2-\www{\ep}, T-\delta_2+\www{\ep}) \le C(M)(D^2 + D^{2\theta}).
$$
Since $\delta_2 - \www{\ep} = (N-1)\www{\ep} = \ep$ and
and $T-\delta_2 + \www{\ep} = T-\ep$, by (3.8), we can complete the proof 
of Theorem 2.
$\blacksquare$

\section{Proof of Theorem 3}

{\bf Section 4.1.  Proof of Theorem 3 (i).}
\\
{\bf First Step.}

We set 
$$
\www{D}_1:= \sum_{j,k=0}^1 
(\Vert \nabla_{x,t}^j \ppp_t^k(\rrr v) \Vert_{L^2(\Gamma \times I)}
+ \Vert \nabla^j\ppp_t^kv\Vert_{L^2(\Gamma\times I)})
$$
and we see
$$
D_1 = \www{D}_1 + \Vert v(\cdot,t_0)\Vert_{H^3(\OOO)}.
$$
\\

We choose $d \in C^2(\ooo{\OOO})$ satisfying (3.1) and (3.2), and we set 
$$
\va(x,t) = e^{\la(d(x) - \beta(t-t_0)^2)}
$$
with sufficiently large $\la>0$ and some $\beta > 0$.  Later we will 
specify $\beta>0$.  By (3.2) we can choose small $\ep_1>0$ such that 
$d=0$ on $\ppp\OOO \setminus \Gamma$, that is,
$$
\va(x,t) \le 1 \quad \mbox{for $(x.t) \in (\ppp\OOO\setminus 
\Gamma)\times I$}, \quad d(x) \ge \ep_1 \quad \mbox{for $x\in \OOO_0$}.
$$
We set 
$$
\mu_1 := e^{\la\ep_1}, \quad 
\mu_2:= \max\{ 1, \, e^{\la(\dddd - \beta\delta^2)}\}.    \eqno{(4.1)}
$$
For $\ep_1>0$, we choose $\beta > 0$ large such that 
$$
\ep_1 > \dddd - \beta \delta^2.                \eqno{(4.2)}
$$
Therefore, (4.1) yields
$$
\left\{ \begin{array}{rl}
& \int_{(\ppp\OOO\setminus \Gamma)\times I}
\sum_{j,k=0}^1 (\vert \nabla_{x,t}^j \ppp_t^k(\rrr v) \vert^2
+ \vert \nabla^j\ppp_t^kv\vert^2) \weight dSdt \le CM_1^2e^{2s\mu_2},\cr\\
& \int_{\Gamma\times I}\sum_{j,k=0}^1 
(\vert \nabla_{x,t}^j \ppp_t^k(\rrr v) \vert^2 
+ \vert \nabla^j\ppp_t^kv\vert^2 \weight) dSdt \le Ce^{Cs}\www{D}_1^2.
\end{array}\right.
                    \eqno{(4.3)}
$$
Moreover, 
$$
C \sum_{\kappa=0}^1\int_{\OOO} \sum_{k=0}^1 
(\vert \nabla (\rrr\ppp_t^kv)(x,t_0+(-1)^{\kappa}\delta)\vert^2
+ \vert \rrr\ppp_t^kv(x,t_0+(-1)^{\kappa}\delta)\vert^2)
e^{2s\va(x,t_0+\delta)} dx
$$
$$
\le CM_1^2e^{2se^{\la(\Vert d\Vert_{C(\ooo{\OOO})} 
- \beta\delta^2)}} \le CM_1^2e^{2s\mu_2}.
                                                         \eqno{(4.4)}
$$

We apply Theorem 1 to $v$ to obtain
$$
\int_Q  (s\vert \nabla v\vert^2 
+ s\vert \nabla(\rrr v)\vert^2 + s^3\vert v\vert^2) 
\weight dxdt 
\le C\int_Q \vert \rrr F \vert^2 \weight dxdt
$$
$$
+ Cs^3M_1^2 e^{2s\mu_2} + Cs^3e^{Cs}\www{D}_1^2                \eqno{(4.5)}
$$
by (4.2) - (4.4).

Setting $v_1:= \ppp_tv$, we have
$$
\left\{\begin{array}{rl}
& \ppp_tv_1(x,t) - \Delta v_1 + (A(x,t)\cdot\nabla)v_1
+ (v_1\cdot \nabla)B(x,t) \cr\\
&+ (\ppp_tA\cdot \nabla)v 
+ (v\cdot\nabla)\ppp_tB + \nabla \ppp_tp = \ppp_tF(x,t),  \cr\\
& \ddd v_1(x,t) = 0, \quad x \in \OOO, \, 0<t<T. 
\end{array}\right.                                       
$$
Noting that 
$$
\rrr((\ppp_tA\cdot \nabla)v) 
= (\ppp_tA\cdot \nabla)\rrr v                        \eqno{(4.6)}
$$
$$
+ \left( \sum_{j=1}^3 (\ppp_2a_j)\ppp_jv_3 
- (\ppp_3a_j)\ppp_jv_2, \, 
\sum_{j=1}^3 (\ppp_3a_j)\ppp_jv_1 - (\ppp_1a_j)\ppp_jv_3, \, 
\sum_{j=1}^3 (\ppp_1a_j)\ppp_jv_2 - (\ppp_2a_j)\ppp_jv_1 \right)^T,
$$
where we set $\ppp_tA:= (a_1, a_2, a_3)^T$, by (4.5) we have
\begin{align*}
& \int_Q (\vert \rrr((\ppp_tA\cdot \nabla)v)\vert^2
+ \vert \rrr((v \cdot \nabla)\ppp_tB)\vert^2) \weight dxdt\\
\le& C\int_Q (\vert \nabla(\rrr v)\vert^2 + \vert \nabla v\vert^2
+ \vert v\vert^2) \weight dxdt
\end{align*}
$$
\le \frac{C}{s}\int_Q \vert \rrr F\vert^2 \weight dxdt
+ Cs^2M_1^2 e^{2s\mu_2} + Cs^2e^{Cs}\www{D}_1^2.             \eqno{(4.7)}
$$
Therefore, we apply Theorem 1 to $v_1$ to obtain
$$
 \int_Q \left( \frac{1}{s}\vert \ppp_t^2\rrr v\vert^2
+ s^3\vert \ppp_t\rrr v\vert^2 
+ s(\vert \nabla (\rrr \ppp_tv)\vert^2 + \vert \nabla \ppp_tv\vert^2)
+ s^3\vert \ppp_tv\vert^2\right) \weight dxdt 
$$
$$
\le C\int_Q \sum_{k=0}^1 \vert \ppp_t^k\rrr F \vert^2 \weight dxdt
+ Cs^3M_1^2 e^{2s\mu_2} + Cs^3e^{Cs}\www{D}_1^2            \eqno{(4.8)}
$$
for all large $s>0$.
\\
{\bf Second Step.}

We start to estimate $\rrr \ppp_tv(x,t_0)$ as follows:
\begin{align*}
& \int_{\OOO} \vert \rrr \ppp_tv(x,t_0)\vert^2 e^{2s\va(x,t_0)}dx
= \int^{t_0}_{t_0-\delta} \left(\ppp_t\int_{\OOO}
\vert \rrr \ppp_tv(x,t)\vert^2 e^{2s\va(x,t)}dx \right) dt\\
+& \int_{\OOO} \vert \rrr \ppp_tv(x,t_0-\delta)\vert^2 
e^{2s\va(x,t_0-\delta)} dx \\
=& \int^{t_0}_{t_0-\delta} \int_{\OOO}
\left(2( (\rrr \ppp_tv)\cdot (\rrr \ppp_t^2v))
+ 2s\ppp_t\va \vert \rrr \ppp_tv\vert^2\right)
e^{2s\va(x,t)} dx dt\\
+& \int_{\OOO} \vert \rrr \ppp_tv(x,t_0-\delta)\vert^2 
e^{2s\va(x,t_0-\delta)} dx \\
\le& C\int_Q (s\vert \rrr \ppp_tv\vert^2
+  \vert \rrr \ppp_tv\vert \vert \rrr \ppp_t^2v\vert)
e^{2s\va(x,t)} dxdt + CM_1^2e^{2s\mu_2}.
\end{align*}
Here we applied also (4.1) and the bound $M_1$.
We have
\begin{align*}
& \vert \rrr \ppp_tv\vert \vert \rrr \ppp_t^2v\vert 
= \left(\frac{1}{s}\vert \rrr \ppp_t^2v\vert\right)
\left( s\vert \rrr \ppp_tv\vert\right)\\
\le & \frac{1}{2}\left( 
\frac{1}{s^2}\vert \rrr \ppp_t^2v\vert^2 + s^2\vert \rrr \ppp_tv\vert^2
\right)
= \frac{1}{2s}\left( 
\frac{1}{s}\vert \rrr \ppp_t^2v\vert^2 + s^3\vert \rrr \ppp_tv\vert^2
\right),
\end{align*}
and so (4.8) implies
\begin{align*}
& \int_Q \vert \rrr \ppp_tv\vert \vert \rrr \ppp_t^2v\vert \weight dxdt
\le & \frac{C}{s}\int_Q \sum_{k=0}^1 \vert \ppp_t^k\rrr F \vert^2 \weight dxdt
+ Cs^2M_1^2 e^{2s\mu_2} + Cs^2e^{Cs}\www{D}_1^2 
\end{align*}
and
$$
\int_Q s\vert \rrr \ppp_tv\vert^2 \weight dxdt 
\le \frac{C}{s^2}\int_Q \sum_{k=0}^1 \vert \ppp_t^k\rrr F \vert^2 \weight dxdt
+ CsM_1^2 e^{2s\mu_2} + Cse^{Cs}\www{D}_1^2.
$$
Hence,
$$
\int_{\OOO} \vert \rrr \ppp_tv(x,t_0)\vert^2 e^{2s\va(x,t_0)}dx
\le \frac{C}{s}\int_Q \sum_{k=0}^1 \vert \ppp_t^k\rrr F \vert^2 \weight dxdt
+ Cs^2M_1^2 e^{2s\mu_2} + Cs^2e^{Cs}D_1^2.        \eqno{(4.9)}
$$

By (1.1) and $\ppp_t\rrr v \in H^1(I;L^2(\OOO)) \subset
C(\ooo{I};L^2(\OOO))$, we have
$$
F(x,t_0) = \ppp_tv(x,t_0) - \Delta v(x,t_0) + (A\cdot \nabla)v(x,t_0)
+ (v(x,t_0)\cdot \nabla)B(x,t_0) + \nabla p(x,t_0),
$$
and 
$$
\rrr F(x,t_0) = \rrr \ppp_tv(x,t_0) + a(x),     \eqno{(4.10)}
$$
where 
$$
a(x) := \rrr(- \Delta v(x,t_0) + (A\cdot \nabla)v(x,t_0)
+ (v(x,t_0)\cdot \nabla)B(x,t_0)).                           \eqno{(4.11)}
$$
Consequently, 
$$
\int_{\OOO} \vert \rrr F(x,t_0)\vert^2 e^{2s\va(x,t_0)} dx
$$
$$ 
\le \frac{C}{s}\int_Q \sum_{k=0}^1 \vert \ppp_t^k\rrr F \vert^2 \weight dxdt
+ Cs^2M_1^2 e^{2s\mu_2} + Cs^2e^{Cs}D_1^2.      
                                              \eqno{(4.12)}
$$
By (1.6) and $\va(x,t) \le \va(x,t_0)$ for $(x,t)\in Q$, we have
$$
\frac{C}{s}\int_Q \sum_{k=0}^1 \vert \rrr \ppp_t^kF(x,t)\vert^2 
e^{2s\va(x,t)} dxdt
\le \frac{C}{s}\int_{\OOO} \vert \rrr F(x,t_0)\vert^2 e^{2s\va(x,t_0)} dx. 
$$
Therefore, 
$$
\int_{\OOO} \vert \rrr F(x,t_0)\vert^2 e^{2s\va(x,t_0)} dx 
\le \frac{C}{s}\int_{\OOO} \vert \rrr F(x,t_0)\vert^2 e^{2s\va(x,t_0)} dx
+ Cs^2M_1^2 e^{2s\mu_2} + Cs^2e^{Cs}D_1^2
$$
for all large $s>0$.

Choosing $s>0$ large, we can absorb the first term on the right-hand side 
into the left-hand side, we obtain
$$
\int_{\OOO} \vert \rrr F(x,t_0)\vert^2 e^{2s\va(x,t_0)} dx 
\le Cs^2M_1^2 e^{2s\mu_2} + Cs^2e^{Cs}D_1^2
$$
for all large $s>s_0$.
By $d(x) \ge \ep_1$ for $x\in \OOO_0$ and $\OOO_0 \subset \OOO$, 
noting $\va(x,t_0) = e^{\la d(x)}$, we see 
$$
e^{2se^{\la\ep_1}}\int_{\OOO_0} \vert \rrr F(x,t_0)\vert^2 dx 
\le \int_{\OOO} \vert \rrr F(x,t_0)\vert^2 e^{2s\va(x,t_0)} dx.
$$
Therefore, 
$$
\int_{\OOO_0} \vert \rrr F(x,t_0)\vert^2 dx 
\le Cs^2M_1^2 e^{2s\mu_2}e^{-2s\mu_1} + Cs^2e^{Cs}D_1^2.
$$ 
By (4.2), we see
$$
\mu_0:= e^{\la\ep_1} - \mu_2 = \mu_1 - \mu_2 > 0,
$$
and so
$$
\int_{\OOO_0} \vert \rrr F(x,t_0)\vert^2 dx 
\le Cs^2M_1^2 e^{-2s\mu_0} + Cs^2e^{Cs}D_1^2
$$ 
for all $s>s_0$.  By the same argument in obtaining (3.11), we complete the 
proof of Theorem 3 (i).
\\
\vspace{0.2cm}
{\bf Section 4.2. Proof of Theorem 3 (ii).}
\\
{\bf First Step.}

We set 
$$
\www{D}_2:= \sum_{j=0}^1 \sum_{k=0}^2 (\Vert \nabla_{x,t}^j\ppp_t^k
(\rrr v)\Vert_{L^2(\Gamma\times I)}
+ \Vert \nabla^j\ppp_t^kv\Vert_{L^2(\Gamma\times I)})
$$
and
$$
D_2 := \www{D}_2 + \Vert v(\cdot,t_0)\Vert_{H^4(\OOO)}, \quad
\www{M}_2 := M_2+m_1.
$$
\\

Now we choose $\va(x,t) = e^{\la(d(x) - \beta(t-t_0)^2)}$, 
where $d\in C^2(\ooo{\OOO})$ satisfies (3.1) and (3.2), and 
$\la>0$ is sufficiently large, and the constant $\beta > 0$ satisfies 
(4.2). 
Then (3.4) and (3.5) hold. 

Therefore, 
$$
\left\{ \begin{array}{rl}
& \int_{(\ppp\OOO\setminus \Gamma)\times I}
\sum_{j=0}^1 \sum_{k=0}^2 (\vert \nabla_{x,t}^j \ppp_t^k(\rrr v) \vert^2
+ \vert \nabla^j\ppp_t^kv\vert^2) \weight dSdt \le CM_2^2e^{2s\mu_2},\cr\\
& \int_{\Gamma\times I}\sum_{j=0}^1 \sum_{k=0}^2 
(\vert \nabla_{x,t}^j \ppp_t^k(\rrr v) \vert^2 
+ \vert \nabla^j\ppp_t^kv\vert^2) \weight dSdt \le Ce^{Cs}\www{D}_2^2.
\end{array}\right.
                    \eqno{(4.13)}
$$

Moreover, 
$$
C \sum_{\kappa=0}^1\int_{\OOO} \sum_{j=0}^1 \sum_{k=0}^2
(\vert \nabla (\rrr\ppp_t^kv)(x,t_0+(-1)^{\kappa}\delta)\vert^2
+ \vert \rrr\ppp_t^kv(x,t_0+(-1)^{\kappa}\delta)\vert^2)
e^{2s\va(x,t_0+\delta)} dx
$$
$$
\le CM_2^2e^{2se^{\la(\Vert d\Vert_{C(\ooo{\OOO})} - \beta\delta^2)}} 
\le CM_2^2e^{2s\mu_2}.
                                                         \eqno{(4.14)}
$$
Moreover, setting $v_2:= \ppp_t^2v$, we have
$$
\left\{\begin{array}{rl}
& \ppp_tv_2(x,t) - \Delta v_2 + (A(x,t)\cdot\nabla)v_2
+ (v_2\cdot \nabla)B(x,t) \cr\\
+ &2(\ppp_tA\cdot \nabla)v_1 + 2(v_1\cdot\nabla)\ppp_tB 
+ (\ppp_t^2A\cdot \nabla)v + (v\cdot\nabla)\ppp_t^2B 
+ \nabla \ppp_t^2p = \ppp_t^2F(x,t),  \cr\\
& \ddd v_2(x,t) = 0, \quad x \in \OOO, \, 0<t<T. 
\end{array}\right.                                       
$$
By (4.13) and (4.14), we apply Theorem 1 to $v_2$ to obtain
\begin{align*}
& \int_Q s\vert \nabla (\rrr \ppp_t^2v)\vert^2 \weight dxdt
\le C\int_Q \vert \rrr \ppp_t^2F\vert^2 \weight dxdt \\
+& C\int_Q  \left\vert \rrr\left\{ 2(\ppp_tA\cdot \nabla)v_1 
+ 2(v_1\cdot\nabla)\ppp_tB 
+ (\ppp_t^2A\cdot \nabla)v 
+ (v\cdot\nabla)\ppp_t^2B\right\} \right\vert^2 \weight dxdt\\
+ & Cs^3M_2^2 e^{2s\mu_2} + Cs^3e^{Cs}\www{D}_2^2
\end{align*}
$$
\le C\int_Q \vert \ppp_t^2\rrr F \vert^2 \weight dxdt
+ Cs^3M_2^2 e^{2s\mu_2} + Cs^3e^{Cs}\www{D}_2^2            \eqno{(4.15)}
$$
for all large $s>0$.

Similarly to (4.6) and (4.7), in terms of (4.8), we have
\begin{align*}
& \int_Q \vert \rrr( 2(\ppp_tA\cdot\nabla) v_1 
+ 2(v_1\cdot \nabla)\ppp_tB)\vert^2 \weight dxdt\\
\le & C\int_Q (\vert \nabla (\rrr v_1)\vert^2
+ \vert \nabla v_1\vert^2 + \vert v_1\vert^2) \weight dxdt\\
= & C\int_Q (\vert \nabla (\rrr \ppp_tv)\vert^2
+ \vert \nabla \ppp_tv\vert^2 + \vert \ppp_tv\vert^2) \weight dxdt
\end{align*}
$$
\le \frac{C}{s}\int_Q \sum_{k=0}^1 \vert \rrr \ppp_t^kF\vert^2 
\weight dxdt
+ Cs^2M_2^2 e^{2s\mu_2} + Cs^2e^{Cs}\www{D}_2^2.        \eqno{(4.16)}
$$
Here we used also $M_1 \le M_2$ and $D_1 \le D_2$.
Similarly we can verify
$$
\int_Q \vert \rrr((\ppp_t^2A\cdot\nabla) v 
+ (v\cdot \nabla)\ppp_t^2B) \vert^2 \weight dxdt
$$
$$
\le \frac{C}{s}\int_Q \vert \rrr F\vert^2 \weight dxdt
+ Cs^2M_2^2 e^{2s\mu_2} + Cs^2e^{Cs}\www{D}_2^2.        \eqno{(4.17)}
$$

Substituting (4.16) and (4.17) into (4.15), we obtain
\begin{align*}
& \int_Q  s\vert \nabla (\rrr\ppp_t^2v)\vert^2 \weight dxdt\\
\le& C\int_Q \sum_{k=0}^2 \vert \rrr \ppp_t^kF\vert^2 \weight dxdt
+ Cs^3M_2^2 e^{2s\mu_2} + Cs^3e^{Cs}\www{D}_2^2. 
\end{align*}
With (4.8), we see
$$
\int_Q (s\vert \nabla(\rrr \ppp_t^2v)\vert^2
+ s\vert \nabla (\rrr \ppp_tv)\vert^2) \weight dxdt
$$
$$
\le C\int_Q \sum_{k=0}^2 \vert \rrr \ppp_t^kF\vert^2 \weight dxdt
+ Cs^3M_2^2 e^{2s\mu_2} + Cs^3e^{Cs}\www{D}_2^2        \eqno{(4.18)}
$$
for all large $s>0$.
\\
{\bf Second Step.}

Now, similarly to Second Step of the proof of Theorem 3 (i), 
using (4.18) and (4.14), we see
\begin{align*}
& \int_{\OOO} \vert \nabla(\rrr \ppp_tv(x,t_0))\vert^2 e^{2s\va(x,t_0)}dx
= \int^{t_0}_{t_0-\delta} \left(\ppp_t\int_{\OOO}
\vert \nabla(\rrr \ppp_tv(x,t))\vert^2 e^{2s\va(x,t)}dx \right) dt\\
+ & \int_{\OOO} \vert \nabla(\rrr \ppp_tv)(x,t_0-\delta)\vert^2
e^{2s\va(x,t_0-\delta)} dx\\
=& \int^{t_0}_{t_0-\delta} \int_{\OOO}
( 2(\nabla(\rrr \ppp_tv)\cdot \nabla(\rrr \ppp_t^2v))
+ 2s\ppp_t\va \vert \nabla(\rrr \ppp_tv)\vert^2)e^{2s\va(x,t)} dx dt\\
+ & \int_{\OOO} \vert \nabla(\rrr \ppp_tv)(x,t_0-\delta)\vert^2
e^{2s\va(x,t_0-\delta)} dx\\
\le& C\int_Q (s\vert \nabla(\rrr \ppp_tv)\vert^2
+  \vert \nabla(\rrr \ppp_tv)\vert \vert \nabla(\rrr \ppp_t^2v)\vert)
e^{2s\va(x,t)} dxdt + CM_2^2e^{2s\mu_2} \\
\le& C\sum_{k=0}^2 \int_Q \vert \rrr \ppp_t^kF\vert^2 \weight dxdt 
+ Cs^3M_2^2 e^{2s\mu_2} + Cs^3e^{Cs}\www{D}_2^2.
\end{align*}
Here $(\nabla(\rrr A)\cdot \nabla(\rrr B))$ means that the sum of the 
products of all the corresponding elements of the two matrices
$\nabla(\rrr A)$ and $\nabla(\rrr B)$.
\\

Applying (1.7) to the first term on the right-hand side, we obtain
$$
\int_{\OOO} \vert \nabla(\rrr \ppp_tv(x,t_0))\vert^2 e^{2s\va(x,t_0)}dx
\le C\int_Q (\vert F(x,t_0)\vert^2 + \vert \nabla F(x,t_0)\vert^2) 
\weight dxdt
$$
$$
+ Cs^3M_2^2 e^{2s\mu_2} + Cs^3e^{Cs}\www{D}_2^2     \eqno{(4.19)}
$$
for all large $s>0$.  

By (1.5), we have
$$
\rrr \rrr F(x,t_0) = -\Delta F(x,t_0) + \nabla (\ddd F(x,t_0))
= -\Delta F(x,t_0), \quad x \in \OOO
$$
and (4.10) yields 
$$
\rrr \rrr F(x,t_0) = \rrr \rrr \ppp_tv(x,t_0)
+ \rrr a(x),
$$
where $a$ is defined by (4.11).
Hence,
$$
-\Delta F(x,t_0) = \rrr \rrr (\ppp_tv(x,t_0)) + \rrr a(x), \quad x\in \OOO.
$$
Therefore, 
$$
\vert \Delta F(x,t_0) \vert 
\le C(\vert \nabla ( \rrr \ppp_tv(x,t_0))\vert 
+ \vert \rrr a(x)\vert, \quad x\in \OOO
$$
and so (4.19) implies 
\begin{align*}
& \int_{\OOO} \vert \Delta F(x,t_0)\vert^2 e^{2s\va(x,t_0)}dx
\le C\int_Q (\vert F(x,t_0)\vert^2 + \vert \nabla F(x,t_0)\vert^2) 
e^{2s\va(x,t)} dxdt \\
+ & Ce^{Cs}(s^3\www{D}_2^2 + \Vert v(\cdot,t_0)\Vert^2_{H^4(\OOO)})
+ Cs^3M_2^2 e^{2s\mu_2}                              
\end{align*}
for all large $s>0$.  Since $\va(x,t) \le \va(x,t_0)$ for $(x,t) \in Q$,
we obtain
$$
\int_{\OOO} \vert \Delta F(x,t_0)\vert^2 e^{2s\va(x,t_0)}dx
\le C\int_{\OOO} (\vert F(x,t_0)\vert^2 + \vert \nabla F(x,t_0)\vert^2) 
e^{2s\va(x,t_0)} dx 
$$
$$
+ Ce^{Cs}s^3 D_2^2 + Cs^3M_2^2 e^{2s\mu_2}            \eqno{(4.20)}             $$
for all large $s>0$
\\
{\bf Third Step.}
\\
Next we apply Lemma 3 shown already in Section 2, which is 
an elliptic Carleman estimate. 
Noting that $F(x,t_0) = \vert \nabla F(x,t_0)\vert = 0$ for 
$x \in \Gamma$, application of Lemma 3 to (4.20) yields
\begin{align*}
&\int_{\OOO} (s\vert \nabla F(x,t_0)\vert^2 + s^3\vert F(x,t_0)\vert^2)
e^{2s\va(x,t_0)}dx \\
\le & C\int_{\OOO} (\vert F(x,t_0)\vert^2 + \vert \nabla F(x,t_0)\vert^2) 
e^{2s\va(x,t_0)} dx 
+ Ce^{Cs}s^3D_2^2 + Cs^3M_2^2e^{2s\mu_2}
\end{align*}
$$
+ Cs^3\int_{\ppp\OOO \setminus \Gamma} 
(\vert \nabla F(x,t_0)\vert^2 + \vert F(x,t_0)\vert^2)e^{2s\va(x,t_0)}dS
                                                                \eqno{(4.21)}
$$
for all large $s>0$.   
Recalling that $\va(x,t) \le 1$ for $(x,t) \in (\ppp\OOO\setminus \Gamma)
\times I$, by $\Vert F(\cdot,t_0)\Vert_{H^2(\OOO)}
\le m_1$ and (4.1), we apply the trace theorem to obtain
$$
\int_{\ppp\OOO \setminus \Gamma} 
(\vert \nabla F(x,t_0)\vert^2 + \vert F(x,t_0)\vert^2)e^{2s\va(x,t_0)}dS
\le m_1^2 \int_{\ppp\OOO\setminus \Gamma} e^{2s} dS
\le Cm_1^2e^{2s} \le Cm_1^2e^{2s\mu_2}.
$$
Choosing $s>0$ large, we can absorb the first term on the right-hand side
of (4.21) into the left-hand side, and so 
$$
\int_{\OOO} s(\vert \nabla F(x,t_0)\vert^2 + \vert F(x,t_0)\vert^2)
e^{2s\va(x,t_0)}dx 
\le Ce^{Cs}s^3D_2^2 + Cs^3\www{M}_2^2e^{2s\mu_2}
$$
for all large $s>0$.
Shrinking the integral domain $\OOO$ on the left-hand side to 
$\OOO_0$ and using (4.1) and
$e^{2s\va(x,t_0)} \ge e^{2se^{\la\ep_1}} = e^{2s\mu_1}$ for $x\in\OOO_0$,
$$
e^{2s\mu_1} \Vert F(\cdot,t_0)\Vert^2_{H^1(\OOO_0)}
\le Ce^{Cs}s^3D_2^2 + Cs^3\www{M}_2^2e^{2s\mu_2}
$$
for all large $s>0$.
By (4.2), we have $\mu_1 > \mu_2$ and so 
$$
\Vert F(\cdot,t_0)\Vert^2_{H^1(\OOO_0)}
\le Ce^{Cs}s^3D_2^2 + Cs^3\www{M}_2^2e^{-2s\mu_0}
$$
for all large $s>0$.  Here we set $\mu_0:= \mu_1 - \mu_2 > 0$.
Therefore, by the same argument as in reaching (3.10), we can complete the 
proof of Theorem 3 (ii).
$\blacksquare$

\section{Concluding Remarks}

{\bf 5.1. Original Navier-Stokes equations}

We mainly consider linearized Navier-Stokes equations.
The original Navier-Stokes equations reads as 
$$
\ppp_tv(x,t) - \Delta v(x,t) + (v\cdot \nabla)v(x,t) + \nabla p 
= F(x,t), \quad x\in \OOO, \, 0<t<T.             \eqno{(5.1)}
$$
Let $v_1$ and $v_2$ satisfy (5.1) with $F_1$ and $F_2$ respectively.
Then, setting $v:= v_1 - v_2$, we have
$$
\ppp_tv(x,t) - \Delta v(x,t) + (v_1\cdot \nabla)v(x,t) 
+ (v\cdot \nabla)v_2(x,t) + \nabla p
= (F_1-F_2)(x,t), \quad x\in \OOO, \, 0<t<T,
$$
which corresponds to (1.1).  Although the existence of $v_1, v_2$ to the 
initial boundary value problems for the Navier-Stokes equations, is
in general, not completely solved with our assumed regularity, for our inverse 
problems we assume the existence of $v_1, v_2$ with such regularity.
Moreover, we see by the proof that $T>0$ can be arbitrarily small, and our
inverse problem requires the existence of the solution $(v,p)$ local in time.
Moreover, we notice that for the inverse problems, we consider 
the solutions only for positive time interval $t> t_0 - \delta > 0$, 
so that there is a possibility that we rely on the smoothing property 
in time of solution to the Navier-Stokes equations, in order to 
gain the necessary regularity of solution for the inverse problem.  

On the other hand, in our inverse source problems, we cannot choose
$t_0=0$, in other words, our problem is not an inverse problem for
any initial boundary value problem.  In particular, we remark
that we do not assume any initial conditions, but we need $v(\cdot,t_0)$ in 
$\OOO$ with some $t_0>0$.  In general, the inverse problem with our 
formulation for parabolic equations as well as the Navier-Stokes 
equations has been a long-standing open problem for the case of
$t_0=0$. 
\\
{\bf 5.2. Data of pressure field $p$}

In this article, we do not use any data of $p$.  Therefore, in Theorem 3 (ii),
we have to assume (1.5): $\ddd F(\cdot,t_0) = 0$ in $\OOO$.
With information of $p$, we need not such assumption for unknown 
sources, but here we do not pursue this direction.
\\
{\bf 5.3. Available Carleman estimates}

As Carleman estimates for parabolic equations including the Navier-Stokes 
equations, we can have two types according to the choices of weight functions
in the following forms.
\begin{itemize}
\item
Regular weight function:
$$
\va(x,t) := e^{\la(d(x) - \beta (t-t_0)^2)}.     \eqno{(5.2)}
$$
\item
Singular weight function:
$$
\va(x,t) = \frac{e^{\la\eta(x)} - e^{2\la\Vert \eta\Vert_{C(\ooo{\OOO})}}}
{h(t)},                       \eqno{(5.3)}
$$
where $\lim_{t\to t_0\pm \delta} h(t) = 0$. 
\end{itemize}
Here $d, \eta \in C^2(\ooo{\OOO})$ are chosen suitably.

The Carleman estimate with the weight (5.3) was proved firstly by 
Imanuvilov \cite{Ima} for a parabolic equation.

As for related inverse problems for the Navier-Stokes equations,
Choulli, Imanuvilov, Puel and Yamamoto \cite{CIPY},
Fan, Di Cristo, Jiang and Nakamura \cite{FDJN}, 
Fan, Jiang and Nakamura \cite{FJN} used Carleman estimates with 
(5.3), while Bellassoued, Imanuvilov and Yamamoto \cite{BIY} and the 
current article rely on (5.2).

Moreover, in order to derive Carleman estimates for the Navier-Stokes equations
with (5.2) or (5.3), we have two options:
\begin{enumerate}
\item
First take $\rrr$ to obtain a parabolic equation in $\rrr v$ and then 
a Poisson equation in $v$. See \cite{CIPY}, \cite{FDJN}, \cite{FJN}
with the weight in form of (5.3).
\item
First take $\ddd$ to obtain a Poisson equation in $p$ and then a
parabolic equation in $v$.  See also \cite{BIY} with the weight 
in form of (5.2).
\end{enumerate}

The present article is based on Option (1) with the weight function (5.2), 
which is a different combination of the weight function and the option.
\\
{\bf 5.4. Improvement of Theorem 3 (ii)}

So far, we prefer to a use of Carleman estimates of a limited category.
Here we apply an $H^{-1}$-Carleman estimate for $\Delta$ and improve 
Theoerm 3 (ii).

Indeed we can relax condition (1.7) as 
$$
\vert \ppp_t^k \rrr F(x,t)\vert \le C(\vert \nabla F(x,t_0)\vert
+ \vert F(x,t_0)\vert), \quad (x,t) \in Q, \, k=0,1.
                                                           \eqno{(5.4)}
$$
Compared with (1.7), we need not the second derivative in $t$ on the 
left-hand side.  Then reducing the time regularity by $1$ in data and 
a priori bound, we can improve the conclusion in Theorem 3 (ii) as
follows.  We recall that $D_1$ is defined by (1.8).
\\

{\bf Proposition 1.}\\
{\it 
There exist constants $C>0$ and $\theta \in (0,1)$, which are dependent on 
$\OOO_0$, $M_1$, $\Gamma$, such that 
$$
\Vert F\Vert_{H^1(\OOO_0)}
\le C( D_1 + D_1^{\theta})                                    \eqno{(5.5)}
$$
for each 
$$
F \in \mathcal{F}_3:= \{ F\in H^1(I;H^2(\OOO));\,
\mbox{supp} F(\cdot,t_0) \subset \OOO, \,  
\mbox{(1.5) and (5.4) are satisfied}\},
$$
provided that (1.9) holds with arbitrarily chosen constant $M_1>0$.
}
\\

We can relax also the condition supp $F(\cdot, t_0) \subset \OOO$, but 
for simplicity we keep this condition in $\mathcal{F}_3$ and here give 
\\
{\bf Sketch of the proof of (5.5).}

Under the assumption, we can follow the argument in the proof of 
Theorem 3 (i) and reach (4.12).
Application of (5.4) to the first term on the right-hand side of 
(4.12) yields
$$
\int_{\OOO} \vert \rrr F(x,t_0)\vert^2 e^{2s\va(x,t_0)} dx 
\le \frac{C}{s}\int_Q (\vert F(x,t_0)\vert^2
+ \vert \nabla F(x,t_0)\vert^2) e^{2s\va(x,t)} dxdt
$$
$$
+ Cs^2M_1^2 e^{2s\mu_2} + Ce^{Cs}s^2D_1^2
                                              \eqno{(5.6)}
$$
for all large $s>0$. 
We estimate the first term on the right-hand side of (5.6) as follows:
\begin{align*}
& \int_Q (\vert F(x,t_0)\vert^2
+ \vert \nabla F(x,t_0)\vert^2) e^{2s\va(x,t)} dxdt\\
= & \int_{\OOO} (\vert F(x,t_0)\vert^2
+ \vert \nabla F(x,t_0)\vert^2) e^{2s\va(x,t_0)} 
\left( \int^{t_0+\delta}_{t_0-\delta} e^{2s(\va(x,t) - \va(x,t_0))}dt
\right) dx.
\end{align*}
Since $d\ge 0$ on $\ooo{\OOO}$, we see 
\begin{align*}
&\int^{t_0+\delta}_{t_0-\delta} e^{2s(\va(x,t) - \va(x,t_0))}dt
= \int^{t_0+\delta}_{t_0-\delta} e^{-2se^{\la d(x)}
(1-e^{-\la\beta (t-t_0)^2})} dt \\
\le & \int^{t_0+\delta}_{t_0-\delta} e^{-2s(1-e^{-\la\beta (t-t_0)^2})} dt
= o(1)
\end{align*}
as $s\to \infty$ by the Lebesgue theorem.  Therefore, 
$$
\int_Q (\vert F(x,t_0)\vert^2 + \vert \nabla F(x,t_0)\vert^2) e^{2s\va(x,t)} 
dxdt
= o(1) \int_{\OOO} (\vert F(x,t_0)\vert^2
+ \vert \nabla F(x,t_0)\vert^2) e^{2s\va(x,t_0)} dx
$$
and so (5.6) yields
$$
s\int_{\OOO} \vert \rrr F(x,t_0)\vert^2 e^{2s\va(x,t_0)} dx
$$
$$
\le o(1) \int_{\OOO} (\vert F(x,t_0)\vert^2
+ \vert \nabla F(x,t_0)\vert^2) e^{2s\va(x,t_0)} dx
+ Cs^3M_1^2 e^{2s\mu_2} + Ce^{Cs}s^3D_1^2         \eqno{(5.7)}
$$
for all large $s>0$.

Now by (1.5), we have
$$
-\Delta F(x,t_0) = \rrr \rrr F(x,t_0), \quad x\in \OOO,  \eqno{(5.8)}
$$
and supp $F(\cdot,t_0) \subset \OOO$.  We show an $H^{-1}$-Carleman 
estimate for $\Delta$ as follows.  
Let $E, \www{E} \subset \R^n$ be bounded domains with smooth boundaries
$\ppp E$ and $\ppp\www{E}$ such that $\ooo{E} \subset \www{E}$, and 
$\omega \subset \www{E} \setminus \ooo{E}$ be an arbitrary non-empty 
open set.  Then in Imanuvilov \cite{Ima} it is proved that there 
exists $\eta \in C^2(\ooo{E})$ satisfying
$$
\eta = 0 \quad \mbox{on $\ppp \www{E}$}, \quad
\eta > 0 \quad \mbox{in $\www{E}$}, \quad
\vert \nabla \eta\vert > 0 \quad \mbox{on $\ooo{\www{E} \setminus \omega}$}.
                                                      \eqno{(5.9)}
$$
We set $\va_0(x) := e^{\la\eta(x)}$ with sufficiently large constant
$\la > 0$.
Then
\\
{\bf Lemma 4.}\\
{\it
There exist constants $C>0$ and $s_0>0$ such that 
$$
\int_E (\vert \nabla w(x)\vert^2 + s^2\vert w(x)\vert^2) e^{2s\va_0(x)} dx
\le C\sum_{j=1}^n s\int_E \vert g_j(x)\vert^2 e^{2s\va_0(x)} dx
$$
for all $s > s_0$ and $w \in H^2_0(E)$ satisfying
$$
\Delta w = \sum_{j=1}^n \ppp_jg_j(x), \quad 
x\in E.
$$
}

This is a Carleman estimate where the right-hand side is estimated in 
$H^{-1}$-space, while Lemma 3 in Section 2 is a Carleman estimate whose
right-hand side is estimated in the $L^2$-space.

Lemma 4 is proved as follows: 
We make the zero extension of $w$ in $E$ to $\www{E} \supset E$.
By $w\in H^2_0(E)$, we see that $w \in H^2_0(\www{E})$ and 
$w=0$ in $\www{E} \setminus E$.  Hence, 
applying Theorem 2.2 in Imanuvilov and Puel \cite{IP}, we have
$$
\int_{\www{E}} (\vert \nabla w(x)\vert^2 + s^2\vert w(x)\vert^2) 
e^{2s\va_0(x)} dx
\le C\sum_{j=1}^n s\int_{\www{E}} \vert g_j(x)\vert^2 e^{2s\va_0(x)} dx
$$
for all large $s>0$, which is Lemma 4.
\\

Here we fix sufficiently large $\la>0$ and so in the Carleman 
estimate we can omit the 
dependence of the constants and the estimate on $\la$.
\\

We apply Lemma 4 by setting $\www{E} = \OOO_1$ and
$E = \OOO$ and $\eta(x) = d(x)$, where $\OOO_1 \supsetneqq \ooo{\OOO}$
satisfies (3.1).  
Moreover $\rrr \rrr F(x,t_0)$ is 
given by a linear combination of $\ppp_j[\rrr F(x,t_0)]_k$,
$1\le j,k\le 3$, where $[\cdot]_k$ denotes the $k$-th component of a vector 
under consideration.  Application of Lemma 4 to (5.8) implies
$$
\int_{\OOO} (\vert \nabla F(x,t_0)\vert^2 + s^2\vert F(x,t_0)\vert^2) 
e^{2s\va(x,t_0)} dx
\le Cs\int_{\OOO} \vert \rrr F(x,t_0)\vert^2 e^{2s\va(x,t_0)} dx
$$
for all large $s>0$.  Substituting (5.7), we obtain
\begin{align*}
& \int_{\OOO} (\vert \nabla F(x,t_0)\vert^2 + s^2\vert F(x,t_0)\vert^2) 
e^{2s\va(x,t_0)} dx\\
\le& o(1) \int_{\OOO} (\vert \nabla F(x,t_0)\vert^2 + \vert F(x,t_0)\vert^2) 
e^{2s\va(x,t_0)} dx
+ Cs^3M_1^2e^{2s\mu_2} + Ce^{Cs}s^3 D_1^2
\end{align*}
for all large $s>0$.

Choosing $s>0$ sufficiently large, we can absorb the first term on the 
right-hand side into the left-hand side, we obtain
$$
\int_{\OOO} (\vert \nabla F(x,t_0)\vert^2 + \vert F(x,t_0)\vert^2) 
e^{2s\va(x,t_0)} dx\\
\le Cs^3M_1^2e^{2s\mu_2} + Ce^{Cs}s^3 D_1^2
$$
for all large $s>0$.
Then, by the same argument as in the final part of the proof of 
Theorem 3 (i), we can complete the proof of (5.5).
$\blacksquare$
\\

In place of Lemma 4, one may apply a Carleman estimate 
(e.g., Vogelsang \cite{Vo}) for a system
$\rrr w = g$ and $\ddd w = 0$, but we do not exploit here.
\section{Appendix. Derivation of Lemma 2 from Lemma 3}

Let $t \in I$ be arbitrarily fixed.  Then, applying Lemma 3 to
$-\Delta r(x,t) = g(x,t)$ in $Q$, we have
$$
\int_{\OOO} (s\vert \nabla r(x,t)\vert^2 + s^3\vert r(x,t)\vert^2) 
e^{2\www{s}\va(x,t_0)} dx
$$
$$
\le C\int_{\OOO} \vert g(x,t)\vert^2 e^{2\www{s}\va(x,t_0)} dx
+ C\www{s}^3\int_{\ppp\OOO} (\vert \nabla r(x,t)\vert^2 
+ \vert r(x,t)\vert^2) e^{2\www{s}\va(x,t_0)} dS          \eqno{(6.1)}
$$
for all $\www{s} > s_1$: some constant.

We set $s_*:= s_1e^{\la\beta \delta^2}$.  Let $s > s_*$.  Then
$$
se^{-\la\beta (t-t_0)^2} \ge se^{-\la\beta \delta^2} > s_1
$$
for all $t \in \ooo{I}$.
Then in (6.1), setting $\www{s}:= se^{-\la\beta(t-t_0)^2}$, we see
\begin{align*}
& \int_{\OOO}(se^{-\la\beta (t-t_0)^2}\vert \nabla r(x,t)\vert^2
+ s^3e^{-3\la\beta(t-t_0)^2 }\vert r(x,t)\vert^2)
e^{2se^{-\la\beta(t-t_0)^2 }\va(x,t_0)}dx\\
\le& C\int_{\OOO} \vert g(x,t)\vert^2 
e^{2se^{-\la\beta (t-t_0)^2}\va(x,t_0)}dx
+ Cs^3e^{-3\la\beta (t-t_0)^2}\int_{\ppp\OOO}
(\vert \nabla r\vert^2 + \vert r\vert^2) 
e^{2se^{-\la\beta(t-t_0)^2}\va(x,t_0)}dS.
\end{align*}
Here, by $e^{-\la\beta(t-t_0)^2} \ge e^{-\la\beta \delta^2}$ for $t\in I$,
we see 
\begin{align*}
& \int_{\OOO}(se^{-\la\beta\delta^2}\vert \nabla r(x,t)\vert^2
+ s^3e^{-3\la\beta\delta^2}\vert r(x,t)\vert^2)
e^{2s\va(x,t)}dx\\
\le& C\int_{\OOO} \vert g(x,t)\vert^2 e^{2s\va(x,t)}dx
+ Cs^3\int_{\ppp\OOO}
(\vert \nabla r\vert^2 + \vert r\vert^2) e^{2s\va(x,t)}dS.
\end{align*}
Multiplying with $e^{3\la\beta \delta^2}$ and
integrating over $t\in I$, we complete the 
derivation of Lemma 2.
$\blacksquare$

\vspace{0.2cm}
\section*{Acknowledgment}
The first author was supported partly by NSF grant DMS 1312900.
The third author was supported by Grant-in-Aid for Scientific Research (S)
15H05740 and Grant-in-Aid (A) 20H00117 of 
Japan Society for the Promotion of Science, 
The National Natural Science Foundation of China
(no. 11771270, 91730303), and the RUDN University 
Strategic Academic Leadership Program.
Most of the article was composed when the first author stayed 
at University of Parma 
as visiting professor and he thanks the university for 
offering an excellent occasion.

\end{document}